\documentclass[letterpaper]{article}
\usepackage{amsmath}
\usepackage{amsthm}
\usepackage{amssymb}
\usepackage[letterpaper]{geometry}
\usepackage{bm}
\usepackage{pstricks}
\usepackage{tikz}
\usepackage{algorithmicx}
\usepackage{algpseudocode}

\newcommand{\nb}{{\mathbf n}}
\newcommand{\pb}{{\mathbf p}}

\newcommand{\qb}{{\mathbf q}}

\newcommand{\RR}{\mathbb R}

\newcommand{\ub}{{\mathbf u}}
\newcommand{\vb}{{\mathbf v}}

\newcommand{\xb}{{\mathbf x}}

\newcommand{\ZZ}{{\mathbb Z}}

\newcommand{\eps}{\varepsilon}
\newcommand{\SLZ}[1]{\text{SL}(#1;\ZZ)}

\theoremstyle{plain}

\theoremstyle{remark}

\title{Cell List Algorithms for Nonequilibrium Molecular Dynamics}
\author{Matthew Dobson\footnote{Department of Mathematics and Statistics, 
University of Massachusetts, 710 N. Pleasant Street, Amherst, MA 01003-9305}
\footnote{Corresponding author: dobson@math.umass.edu} \and Ian Fox$^*$ \and Alexandra Saracino$^*$}

\begin{document}
\maketitle
\begin{abstract}
We present two modifications of the standard cell list algorithm
for nonequilibrium molecular dynamics simulations of homogeneous, linear flows.
When such a flow is modeled with periodic boundary conditions, the simulation
box deforms with the flow, and recent progress has been made developing
boundary conditions suitable for general 3D flows of
this type.  For the typical case of short-ranged, pairwise interactions, the
cell list algorithm reduces computational complexity of the force computation
from O($N^2$) to O($N$), where $N$ is the total number of particles in the
simulation box.  The new versions of the cell list algorithm handle the
dynamic, deforming simulation geometry.  We include a comparison of the
complexity and efficiency of the two proposed modifications of the standard
algorithm.
\end{abstract}

\section{Introduction}
Recent developments in the simulation of nonequilibrium molecular dynamics for homogeneous flow have greatly increased the types of flows that can be simulated for long simulation times~\cite{hunt13, dobs14}.  In a simulation of homogeneous, linear background flow with periodic boundary conditions, the simulation box deforms with the background flow.  If care is not taken in the alignment of the simulation box, this can  lead to a breakdown of the numerics.  However, boundary conditions have been devised for shear flow~\cite{Lees}, planar elongational flow~\cite{todd98,bara99,todd99} based on the original study of Kraynik and Reinelt~\cite{Kraynik}, and subsequently for most 3D flows in~\cite{hunt13, dobs14}.  

The deforming geometry of the simulation box presents new challenges when developing efficient algorithms.  Many techniques have been developed to increase simulation speed and scalability for the equilibrium case with a static simulation box.  One prominent example is the cell list algorithm, which reduces the computational complexity of the simulation from O($N^2$) to O($N$) by exploiting the short range interaction of the particles to eliminate unnecessary computation.  In this paper we present two new versions of the cell list algorithm generalized for simulations using the generalized Kraynik-Reinelt boundary conditions.  Previous work by Matin et al.~\cite{Matin} has produced cell list algorithms for planar shear and elongational flows, which used precomputed cell neighborhoods  based on the specific form of the flow. The versions described here allow for robust, dynamic handling of general 3D flows including planar, uniaxial, and biaxial flows. 

In Section~\ref{sec:pbc}, we described periodic boundary conditions for equilibrium and nonequilibrium flows.  In Section~\ref{sec:eqcell}, we review the cell list algorithm for equilibrium molecular dynamics, before introducing the new modifications for nonequilibrium flows in Section~\ref{sec:neqcell}.  Two variants are described, the dynamic size cell list and the dynamic offset cell list, and the efficiency of these two algorithms is compared in Section~\ref{sec:eff}.

\section{Periodic Boundary Conditions}
\label{sec:pbc}

We review the formulation of periodic boundary conditions for equilibrium and nonequilibrium molecular dynamics in order to fix notation and describe the setting for the cell list algorithms presented here.
We first consider an equilibrium molecular dynamics 
simulation with periodic boundary conditions.  The simulation geometry is defined 
by three basis vectors written together in matrix form:
\begin{equation}L = [\vb_1 \  \vb_2 \  \vb_3 ] \in \RR^{3 \times 3}.\end{equation}  These vectors define a lattice in $\RR^3$
whose points are given by $\{L \nb = n_1 \vb_1 + n_2 \vb_2 + n_3 \vb_3 \ | \ \nb \in \ZZ^3\}.$  
The simulation box is the unit cell of the lattice, 
\begin{equation}\Omega = \{ \lambda_1 \vb_1 + \lambda_2 \vb_2 + \lambda_3 \vb_3 \ | \ 0 \leq \lambda_1,\lambda_2,\lambda_3 < 1  \}.\end{equation}
We track the position of each particle within the unit cell, and a particle with position $\qb$ has an infinite number of images, whose positions are at $\qb + L \nb,$ for all $\nb \in \ZZ^3.$  
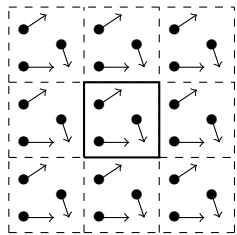
\begin{figure}[ht]
\centering
\begin{tikzpicture}
\draw [dashed] (1,0) -- (3,0);
\draw [dashed] (2,1) -- (3,1);
\draw [dashed] (0,2) -- (3,2);
\draw [dashed] (0,3) -- (3,3);
\draw [dashed] (0,2) -- (0,3);
\draw [dashed] (1,0) -- (1,3);
\draw [dashed]  (2,0) -- (2,3);
\draw [dashed] (3,0) -- (3,3);
\draw [dashed] (0,0) -- (1,0);
\draw [dashed] (0,0) -- (0,1);
\draw [dashed] (0,1) -- (1,1);
\draw [dashed] (0,1) -- (0,2);
\draw[thick] (1,1) -- (2,1) -- (2,2) -- (1,2) -- (1,1);
\node at (0.2, 2.2) (node) {$\bullet$};
\node at (0.2, 2.7) (node) {$\bullet$};
\node at (0.7, 2.5) (node) {$\bullet$};
\draw[->] (0.2, 2.2) -- (0.6, 2.2);
\draw[->] (0.2, 2.7) -- (0.5, 2.9);
\draw[->] (0.7, 2.5) -- (0.8, 2.2);

\begin{scope}[yshift = -2cm]
\node at (0.2, 2.2) (node) {$\bullet$};
\node at (0.2, 2.7) (node) {$\bullet$};
\node at (0.7, 2.5) (node) {$\bullet$};
\draw[->] (0.2, 2.2) -- (0.6, 2.2);
\draw[->] (0.2, 2.7) -- (0.5, 2.9);
\draw[->] (0.7, 2.5) -- (0.8, 2.2);
\end{scope}

\begin{scope} [yshift = -1cm]
\node at (0.2, 2.2) (node) {$\bullet$};
\node at (0.2, 2.7) (node) {$\bullet$};
\node at (0.7, 2.5) (node) {$\bullet$};
\draw[->] (0.2, 2.2) -- (0.6, 2.2);
\draw[->] (0.2, 2.7) -- (0.5, 2.9);
\draw[->] (0.7, 2.5) -- (0.8, 2.2);
\end{scope}

\begin{scope}[xshift = 1cm]
\node at (0.2, 2.2) (node) {$\bullet$};
\node at (0.2, 2.7) (node) {$\bullet$};
\node at (0.7, 2.5) (node) {$\bullet$};
\draw[->] (0.2, 2.2) -- (0.6, 2.2);
\draw[->] (0.2, 2.7) -- (0.5, 2.9);
\draw[->] (0.7, 2.5) -- (0.8, 2.2);
\end{scope}

\begin{scope}[xshift = 1cm, yshift = -1cm]
\node at (0.2, 2.2) (node) {$\bullet$};
\node at (0.2, 2.7) (node) {$\bullet$};
\node at (0.7, 2.5) (node) {$\bullet$};
\draw[->] (0.2, 2.2) -- (0.6, 2.2);
\draw[->] (0.2, 2.7) -- (0.5, 2.9);
\draw[->] (0.7, 2.5) -- (0.8, 2.2);
\end{scope}

\begin{scope}[xshift = 1cm, yshift = -2cm]
\node at (0.2, 2.2) (node) {$\bullet$};
\node at (0.2, 2.7) (node) {$\bullet$};
\node at (0.7, 2.5) (node) {$\bullet$};
\draw[->] (0.2, 2.2) -- (0.6, 2.2);
\draw[->] (0.2, 2.7) -- (0.5, 2.9);
\draw[->] (0.7, 2.5) -- (0.8, 2.2);
\end{scope}

\begin{scope}[xshift = 2cm]
\node at (0.2, 2.2) (node) {$\bullet$};
\node at (0.2, 2.7) (node) {$\bullet$};
\node at (0.7, 2.5) (node) {$\bullet$};
\draw[->] (0.2, 2.2) -- (0.6, 2.2);
\draw[->] (0.2, 2.7) -- (0.5, 2.9);
\draw[->] (0.7, 2.5) -- (0.8, 2.2);
\end{scope}

\begin{scope}[xshift = 2cm, yshift = -1cm]
\node at (0.2, 2.2) (node) {$\bullet$};
\node at (0.2, 2.7) (node) {$\bullet$};
\node at (0.7, 2.5) (node) {$\bullet$};
\draw[->] (0.2, 2.2) -- (0.6, 2.2);
\draw[->] (0.2, 2.7) -- (0.5, 2.9);
\draw[->] (0.7, 2.5) -- (0.8, 2.2);
\end{scope}

\begin{scope}[xshift = 2cm, yshift = -2cm]
\node at (0.2, 2.2) (node) {$\bullet$};
\node at (0.2, 2.7) (node) {$\bullet$};
\node at (0.7, 2.5) (node) {$\bullet$};
\draw[->] (0.2, 2.2) -- (0.6, 2.2);
\draw[->] (0.2, 2.7) -- (0.5, 2.9);
\draw[->] (0.7, 2.5) -- (0.8, 2.2);
\end{scope}
\end{tikzpicture}
\caption{\label{fig:eqpbc}Periodic boundary conditions at equilibrium. The center square shows a simulation box with three particles. Image particles all experience identical forces and move in sync for all time.}
\end{figure}
In an equilibrium molecular dynamics simulation, each of the image particles has the same velocity, so that for all times the images move in sync with one another.  Thus, a particle with phase coordinates $(\qb, \pb)$ has images with coordinates $(\qb + L \nb, \pb)$ for all $\nb \in \ZZ^3$.

\subsection{Nonequilibrium flow and the deforming simulation box}
In a nonequilibrium simulation of steady, homogeneous flow, there is a background flow represented by the flow matrix $A \in \RR^{3 \times 3}.$  This represents a linear flow field where at spatial position $\xb \in \RR^3,$ the macroscopic velocity equals $\ub(\xb) = A \xb.$  We choose periodic boundary conditions where the image particle velocities are consistent with the background flow. The lattice vectors are now necessarily time dependent and are denoted:
\begin{equation}
L_t = [\vb_{1,t} \  \vb_{2,t} \ \vb_{3,t} ],
\end{equation}
and in this case a particle with phase coordinates $(\qb, \pb)$ has images with coordinates 
$(\qb + L_t \nb, \pb + A L_t \nb)$ for $\nb \in \ZZ^3.$

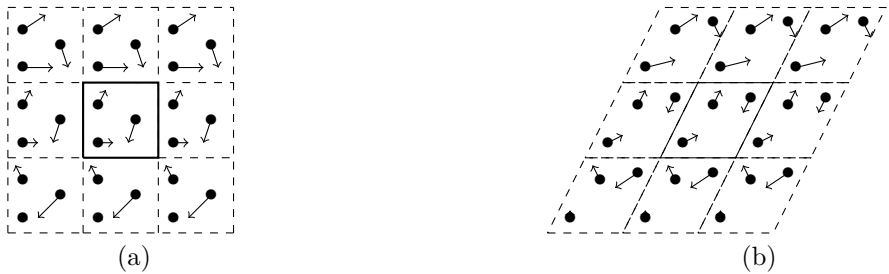
\begin{figure}[ht]
\centerline{
\centering
\begin{tikzpicture}
\draw [dashed] (1,0) -- (3,0);
\draw [dashed] (2,1) -- (3,1);
\draw [dashed] (0,2) -- (3,2);
\draw [dashed] (0,3) -- (3,3);
\draw [dashed] (0,2) -- (0,3);
\draw [dashed] (1,0) -- (1,3);
\draw [dashed]  (2,0) -- (2,3);
\draw [dashed] (3,0) -- (3,3);
\draw [dashed] (0,0) -- (1,0);
\draw [dashed] (0,0) -- (0,1);
\draw [dashed] (0,1) -- (1,1);
\draw [dashed] (0,1) -- (0,2);
\draw[thick] (1,1) -- (2,1) -- (2,2) -- (1,2) -- (1,1);
\node at (0.2, 2.2) (node) {$\bullet$};
\node at (0.2, 2.7) (node) {$\bullet$};
\node at (0.7, 2.5) (node) {$\bullet$};
\draw[->] (0.2, 2.2) -- (0.6, 2.2);
\draw[->] (0.2, 2.7) -- (0.5, 2.9);
\draw[->] (0.7, 2.5) -- (0.8, 2.2);

\begin{scope}[yshift = -2cm]
\node at (0.2, 2.2) (node) {$\bullet$};
\node at (0.2, 2.7) (node) {$\bullet$};
\node at (0.7, 2.5) (node) {$\bullet$};
\draw[->] (0.2, 2.2) -- (0.2, 2.2);
\draw[->] (0.2, 2.7) -- (0.1, 2.9);
\draw[->] (0.7, 2.5) -- (0.4, 2.2);
\end{scope}

\begin{scope} [yshift = -1cm]
\node at (0.2, 2.2) (node) {$\bullet$};
\node at (0.2, 2.7) (node) {$\bullet$};
\node at (0.7, 2.5) (node) {$\bullet$};
\draw[->] (0.2, 2.2) -- (0.4, 2.2);
\draw[->] (0.2, 2.7) -- (0.3, 2.9);
\draw[->] (0.7, 2.5) -- (0.6, 2.2);
\end{scope}

\begin{scope}[xshift = 1cm]
\node at (0.2, 2.2) (node) {$\bullet$};
\node at (0.2, 2.7) (node) {$\bullet$};
\node at (0.7, 2.5) (node) {$\bullet$};
\draw[->] (0.2, 2.2) -- (0.6, 2.2);
\draw[->] (0.2, 2.7) -- (0.5, 2.9);
\draw[->] (0.7, 2.5) -- (0.8, 2.2);
\end{scope}

\begin{scope}[xshift = 1cm, yshift = -1cm]
\node at (0.2, 2.2) (node) {$\bullet$};
\node at (0.2, 2.7) (node) {$\bullet$};
\node at (0.7, 2.5) (node) {$\bullet$};
\draw[->] (0.2, 2.2) -- (0.4, 2.2);
\draw[->] (0.2, 2.7) -- (0.3, 2.9);
\draw[->] (0.7, 2.5) -- (0.6, 2.2);
\end{scope}

\begin{scope}[xshift = 1cm, yshift = -2cm]
\node at (0.2, 2.2) (node) {$\bullet$};
\node at (0.2, 2.7) (node) {$\bullet$};
\node at (0.7, 2.5) (node) {$\bullet$};
\draw[->] (0.2, 2.2) -- (0.2, 2.2);
\draw[->] (0.2, 2.7) -- (0.1, 2.9);
\draw[->] (0.7, 2.5) -- (0.4, 2.2);
\end{scope}

\begin{scope}[xshift = 2cm]
\node at (0.2, 2.2) (node) {$\bullet$};
\node at (0.2, 2.7) (node) {$\bullet$};
\node at (0.7, 2.5) (node) {$\bullet$};
\draw[->] (0.2, 2.2) -- (0.6, 2.2);
\draw[->] (0.2, 2.7) -- (0.5, 2.9);
\draw[->] (0.7, 2.5) -- (0.8, 2.2);
\end{scope}

\begin{scope}[xshift = 2cm, yshift = -1cm]
\node at (0.2, 2.2) (node) {$\bullet$};
\node at (0.2, 2.7) (node) {$\bullet$};
\node at (0.7, 2.5) (node) {$\bullet$};
\draw[->] (0.2, 2.2) -- (0.4, 2.2);
\draw[->] (0.2, 2.7) -- (0.3, 2.9);
\draw[->] (0.7, 2.5) -- (0.6, 2.2);
\end{scope}

\begin{scope}[xshift = 2cm, yshift = -2cm]
\node at (0.2, 2.2) (node) {$\bullet$};
\node at (0.2, 2.7) (node) {$\bullet$};
\node at (0.7, 2.5) (node) {$\bullet$};
\draw[->] (0.2, 2.2) -- (0.2, 2.2);
\draw[->] (0.2, 2.7) -- (0.1, 2.9);
\draw[->] (0.7, 2.5) -- (0.4, 2.2);
\end{scope}
\end{tikzpicture} \hspace{1.5in}
\centering
\begin{tikzpicture}
\begin{scope}
\draw (0,0)--(1,0)--(1.5,1)--(.5,1)--(0,0);
\node at (.3, .2) (node) {$\bullet$};
\node at (.7, .7) (node) {$\bullet$};
\node at (1.2, .8) (node) {$\bullet$};
\draw[->] (.3, .2) -- (.5, .3);
\draw[->] (.7, .7) -- (.8, .9);
\draw[->] (1.2, .8) -- (1.1, .6);
\end{scope}

\begin{scope}[xshift = -1cm]
\draw [dashed](0,0)--(1,0)--(1.5,1)--(.5,1)--(0,0);
\node at (.3, .2) (node) {$\bullet$};
\node at (.7, .7) (node) {$\bullet$};
\node at (1.2, .8) (node) {$\bullet$};
\draw[->] (.3, .2) -- (.5, .3);
\draw[->] (.7, .7) -- (.8, .9);
\draw[->] (1.2, .8) -- (1.1, .6);
\end{scope}

\begin{scope}[xshift = 1cm]
\draw [dashed](0,0)--(1,0)--(1.5,1)--(.5,1)--(0,0);
\node at (.3, .2) (node) {$\bullet$};
\node at (.7, .7) (node) {$\bullet$};
\node at (1.2, .8) (node) {$\bullet$};
\draw[->] (.3, .2) -- (.5, .3);
\draw[->] (.7, .7) -- (.8, .9);
\draw[->] (1.2, .8) -- (1.1, .6);
\end{scope}

\begin{scope}[xshift = -.5cm, yshift = -1cm]
\draw [dashed](0,0)--(1,0)--(1.5,1)--(.5,1)--(0,0);
\node at (.3, .2) (node) {$\bullet$};
\node at (.7, .7) (node) {$\bullet$};
\node at (1.2, .8) (node) {$\bullet$};
\draw[->] (.3, .2) -- (.3, .3);
\draw[->] (.7, .7) -- (.6, .9);
\draw[->] (1.2, .8) -- (.9, .6);
\end{scope}

\begin{scope}[xshift = .5cm, yshift = 1cm]
\draw [dashed](0,0)--(1,0)--(1.5,1)--(.5,1)--(0,0);
\node at (.3, .2) (node) {$\bullet$};
\node at (.7, .7) (node) {$\bullet$};
\node at (1.2, .8) (node) {$\bullet$};
\draw[->] (.3, .2) -- (.7, .3);
\draw[->] (.7, .7) -- (1, .9);
\draw[->] (1.2, .8) -- (1.3, .6);
\end{scope}

\begin{scope}[xshift = -1.5cm, yshift = -1cm]
\draw [dashed](0,0)--(1,0)--(1.5,1)--(.5,1)--(0,0);
\node at (.3, .2) (node) {$\bullet$};
\node at (.7, .7) (node) {$\bullet$};
\node at (1.2, .8) (node) {$\bullet$};
\draw[->] (.3, .2) -- (.3, .3);
\draw[->] (.7, .7) -- (.6, .9);
\draw[->] (1.2, .8) -- (.9, .6);
\end{scope}

\begin{scope}[xshift = .5cm, yshift = -1cm]
\draw [dashed](0,0)--(1,0)--(1.5,1)--(.5,1)--(0,0);
\node at (.3, .2) (node) {$\bullet$};
\node at (.7, .7) (node) {$\bullet$};
\node at (1.2, .8) (node) {$\bullet$};
\draw[->] (.3, .2) -- (.3, .3);
\draw[->] (.7, .7) -- (.6, .9);
\draw[->] (1.2, .8) -- (.9, .6);
\end{scope}

\begin{scope}[xshift = -.5cm, yshift = 1cm]
\draw [dashed](0,0)--(1,0)--(1.5,1)--(.5,1)--(0,0);
\node at (.3, .2) (node) {$\bullet$};
\node at (.7, .7) (node) {$\bullet$};
\node at (1.2, .8) (node) {$\bullet$};
\draw[->] (.3, .2) -- (.7, .3);
\draw[->] (.7, .7) -- (1, .9);
\draw[->] (1.2, .8) -- (1.3, .6);
\end{scope}

\begin{scope}[xshift = 1.5cm, yshift = 1cm]
\draw [dashed](0,0)--(1,0)--(1.5,1)--(.5,1)--(0,0);
\node at (.3, .2) (node) {$\bullet$};
\node at (.7, .7) (node) {$\bullet$};
\node at (1.2, .8) (node) {$\bullet$};
\draw[->] (.3, .2) -- (.7, .3);
\draw[->] (.7, .7) -- (1, .9);
\draw[->] (1.2, .8) -- (1.3, .6);
\end{scope}

\end{tikzpicture}}
\centerline{(a) \hspace{3in} (b)}
\caption{\label{fig:defparr}Periodic boundary conditions under a shear flow.  In (a), the initial configuration is shown, with a square simulation box.  Note that the image particle have differences in the horizontal component of the velocity depending on the vertical position, consistent with the background shear flow.  In b), the simulation box and all replicas have deformed along with the flow.  Again, replica velocities are position dependent, consistent with the background flow.}
\end{figure}

The particles in the simulation box obey the velocity equation $\frac{d \qb}{dt} = m^{-1} \pb,$ and all images likewise satisfy
\begin{equation}\frac{d(\qb + L_t \nb)}{dt} = m^{-1} \pb + A L_t \nb. \end{equation}
Therefore, the lattice vectors move with the flow, 
\begin{equation}\frac{d}{dt}L_t = A L_t \text{, which has solution } L_t = e^{At}L_0.\end{equation}
Thus, the simulation box deforms with the flow, and this deformation can lead to the basis vectors growing 
linearly, as in the case of shear flow, or exponentially, as in the case of elongational flow.

\subsection{Controlling Simulation Box Deformation}
Throughout a simulation it is important that the simulation box does not become too deformed.  Very long vectors can lead to round-off errors and difficulties in computing periodic distances.  More problematically, uncontrolled elongational flow can cause particles to approach their own images, causing numerical instability. Several techniques have been developed to deal with the deformation of the simulation box, and in this work we focus on the class of techniques involving lattice remappings.  

Let $\SLZ{3}$ denote the set of all integer three by three matrices with determinant one. For $M \in \SLZ{3}$ the lattice generated by $L_t$ is equivalent to the lattice generated by $L_t M,$ since $M$ is a one-to-one map of $\ZZ^3$ onto $\ZZ^3.$  Thus we can apply $M$ to choose new lattice basis vectors without changing system dynamics. Using these remappings we choose an initial lattice $L_0$ and at appropriate times during the simulation, remap the lattice basis $\widetilde{L}_t = L_t M,$ where the unit cell of $\widetilde{L}_t$ is less deformed than that of $L_t.$

\begin{figure}[ht]
\centering
\begin{tikzpicture}
\begin{scope}[yshift = 2.5cm]
\fill (-1,0) circle (.1);
\fill (-1,1) circle (.1);
\fill (0,0) circle (.1);
\fill (0,1) circle (.1);
\fill (1,0) circle (.1);
\fill (1,1) circle (.1);
\fill (2,0) circle (.1);
\fill (2,1) circle (.1);
\draw (0,0)--(0,1)--(1,1)--(1,0)--(0,0);
\end{scope}
\begin{scope}[xshift = -2.5cm]
\fill (-1,0) circle (.1);
\fill (-1,1) circle (.1);
\fill (0,0) circle (.1);
\fill (0,1) circle (.1);
\fill (1,0) circle (.1);
\fill (1,1) circle (.1);
\fill (2,0) circle (.1);
\fill (2,1) circle (.1);
\draw (0,0)--(.5,1)--(1.5,1)--(1,0)--(0,0);
\end{scope}
\begin{scope}
\draw[->, very thick] (-.5, 2) -- (-1, 1.5);
\draw[->, very thick] (0, .5) -- (1, .5);
\end{scope}
\begin{scope}[xshift = 2.5cm]
\fill (-1,0) circle (.1);
\fill (-1,1) circle (.1);
\fill (0,0) circle (.1);
\fill (0,1) circle (.1);
\fill (1,0) circle (.1);
\fill (1,1) circle (.1);
\fill (2,0) circle (.1);
\fill (2,1) circle (.1);
\draw (0,0)--(-.5,1)--(.5,1)--(1,0)--(0,0);
\end{scope}
\end{tikzpicture}
\caption{\label{fig:shearres} Resetting the simulation box under shear flow. The simulation box begins undeformed and undergoes shear flow. After it has been sheared a certain amount to the right we can redefine the simulation box as sheared to the left without changing the resulting lattice by applying a lattice automorphism. By repeating this process we ensure the simulation box never becomes too deformed.}
\end{figure}
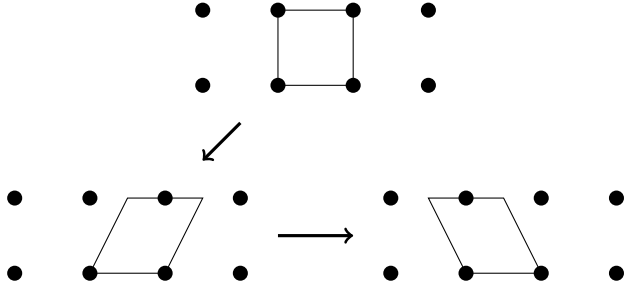

For planar shear flows, with background flow 
\begin{equation}
A = \left[
\begin{array}{rrr}
0 & \eps & 0 \\
0 & 0 & 0 \\
0 & 0 & 0 
\end{array}
\right], 
\end{equation} we can use Lees-Edwards boundary conditions~\cite{Lees, Evans} to perform the remapping.  The initial lattice vectors are parallel to the coordinate axes and satisfy for all time 
\begin{equation}
L_t = \left[
\begin{array}{rrr}
a & a \eps t & 0 \\
0 & a & 0 \\
0 & 0 & a 
\end{array}
\right],
\end{equation} where $a$ is the initial length of the three lattice vectors.  At time $t^* = \frac{1}{\eps},$ the simulation box has been sheared by 45$^\circ,$ and it produces the same lattice as the undeformed box. The simulation box can then be redefined by applying an automorphism to its undeformed state, taking $\widetilde{L_t} = L_t M$ where 
\begin{equation}
\label{eq:Mshear}
M = \left[
\begin{array}{rrr}
1 & -1 & 0 \\
0 & 1 & 0 \\
0 & 0 & 1 
\end{array}
\right] 
\end{equation}
A smaller maximum angle of deformation can be achieved by instead applying $M$ in~\eqref{eq:Mshear} when the simulation box has been sheared one half of a lattice spacing, as shown in Figure~\ref{fig:shearres}. This will produce a maximum angle of deformation equal to $90^{\circ}-tan^{-1}(2) \approx 26.57^{\circ}.$ 

Kraynik-Reinelt (KR) boundary conditions can be used for planar elongational flows~\cite{Kraynik, todd98, bara99}.  The original lattice basis for the simulation box is chosen such that after a finite amount of time the deformed simulation box lattice is mapped to its original configuration.
In other words we choose an $L_0$ such that at time $t_*,$
\begin{equation}\label{eq:kr}L_{t_*} = e^{At_*}L_0 = L_0M\end{equation} 
This is possible for general planar flows~\cite{Kraynik, hunt10}, but has been shown to be impossible for certain 3D flows, including uniaxial and biaxial stretching flows~\cite{Kraynik}. Note the choice of $M$ is made independent of $A$.

To generalize the KR boundary conditions to linear, incompressible flows in three dimensions we use a pair of 
automorphisms that are used to keep the deformation bounded to a limited stretch space.  While the deformation of the simulation box is no longer time periodic, the boundary conditions ensure that $\widetilde{L}_t$ is related to $L_0$ by a bounded stretch.  For example, when $A$ is diagonal, we have
\begin{equation}
\widetilde{L}_t = \exp \left(  \begin{bmatrix}
\eps_{t,1} & 0 & 0 \\
0 & \eps_{t,2} & 0 \\
0 & 0 & \eps_{t,3}
\end{bmatrix} \right) \  L_0,
\end{equation} where the coordinate stretches satisfy $\eps_{t,1} + \eps_{t,2} + \eps_{t,3} = 0,$ and the vector $[\eps_{t,1}, \eps_{t,2}, \eps_{t,3}]$ lives in a bounded parallelogram in $\RR^3.$  These generalized KR boundary conditions are developed more completely in \cite{dobs14}.

\section{Equilibrium Cell List Algorithm}
\label{sec:eqcell}
Force interactions in molecular dynamics simulations are often short-ranged and are implemented with a cutoff radius past which the forces are identically zero.  The naive method of computing pairwise force interactions in the simulation is to loop over all pairs of particles, which leads to $O(N^2)$ time complexity.  The cell list algorithm takes advantage of the short-ranged interactions and reduces overall complexity to $O(N)$ by ignoring the majority of non-contributing pairwise interactions.

In the cell list algorithm the simulation box is divided into many small boxes called cells~\cite{fren01}. Each side length of a cell is at least as large as the radius of interaction. At each simulation step before force interactions are computed each particle is assigned to a cell based off of its position in the simulation box. For a given cell, particles will only have non-zero interactions with particles in the same cell and those in the cells surrounding it, so rather than looping over each particle pair in the simulation box, one only loops over the particle pairs in these 27 cells.  A cell and its 26 adjacent cells are together called a cell neighborhood.  Figure~\ref{fig:eqpbc} shows a small cell list, where for convenience this and all following diagrams are shown in 2D while the algorithms presented here are formulated and tested in 3D.

\begin{figure}[ht]
\centering
\begin{tikzpicture}

\draw[step=.5 cm,gray]
(0,0) grid (2,2);
\foreach \x in {.1, .4, .7,    1.3, 1.6, 1.9}
\foreach \y in {.1, .4, .7, 1, 1.3, 1.6, 1.9}
{
\fill (\x + rand/20, \y + rand/20) circle (.05);
}

\foreach \x in {1}
\foreach \y in {.1, .4, .7, 1, 1.6, 1.9}
{
\fill (\x + rand/20, \y + rand/20) circle (.05);
}

\fill (1, 1.3) circle (.05);
\draw (1, 1.3) circle (.5);

\end{tikzpicture}
\caption{\label{fig:eqclpr}An equilibrium cell list. The simulation box has been divided into sixteen equal cells. The large circle in the middle shows the interaction range around the center particle.  Any particles outside that range have zero force interaction with the center particle.}
\end{figure}
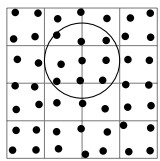

While this algorithm uses $O(N)$ time for dividing and assigning particles at each step in the simulation it reduces the number of force checks required for each particle. Notably the number of checks per particle becomes independent of the number of particles in the system. This reduces the time complexity of interparticle force calculation to $O(N)$.

\section{Nonequilibrium Cell List Algorithm}
\label{sec:neqcell}
The traditional cell list algorithm needs to be modified when simulating nonequilibrium flows to account for the deformation of the simulation box.  We describe two variants of the standard cell list algorithm to deal with this problem.

In the first method we choose cells that are congruent to the simulation box, so that the box is neatly divided into cells. In our approach we seek to keep a standard cell neighborhood of 27 cells, so we choose the cells so that all three cell heights is at least as large as the interaction range.  The shape of the individual cells evolves with the flow, but the total number of cells, and therefore their size, changes in order to satisfy the minimum height requirement.  We call this modified algorithm the dynamic size cell list algorithm. 

In the second method cells remain rectangular boxes. This creates fractional cells along the faces of the simulation box which are combined into complete cells.  The cells do not neatly align across all boundaries, so that in this method offsets between periodic replicas of the simulation box must be considered and the cell neighborhood expanded as necessary across the boundary.  We call this modified algorithm the dynamic offset cell list algorithm.

Both variations of the cell list algorithm we created were developed to handle general three dimensional linear flows. Thus our algorithm handles shear and planar elongational flows in the same manner.  In Section~\ref{sec:eff} we will compare the computational efficiency, overhead, and algorithmic complexity of these variations. 

Previous work has produced cell list algorithm variations for handling shear flow and planar elongational flow~\cite{Matin}. The method developed by Matin et al. allowed cells to deform with the background flow and precomputed an expanded cell neighborhood so that all nonzero particle interactions were included even when cells were maximally deformed. The method made use of a precomputed neighbor list for handling planar flows. While our method uses similar ideas, our use of a dynamically computed neighborhood increases the efficiency of the code and allows for more general treatment of more general flows. 

\subsection{Dynamic Size}

\begin{figure}[ht]
\centering
\begin{tikzpicture}

\begin{scope}[xslant = .5]
\draw[xstep=.667 cm,ystep=.5 cm,black]
(0,0) grid (2.01,2);
\end{scope}

\begin{scope}
\clip (0,0)--(2,0) -- (3,2) -- (1, 2) -- (0, 0);
\foreach \x in {.1, .4, .7, 1, 1.3, 1.9, 2.2, 2.5, 2.8, 3.1, 3.4}
\foreach \y in {.1, .4, .7, 1, 1.3, 1.6, 1.9}
{
\fill (\x + rand/20, \y + rand/20) circle (.05);
}

\foreach \x in {1.6}
\foreach \y in {.1, .4, .7, 1, 1.6, 1.9}
{
\fill (\x + rand/20, \y + rand/20) circle (.05);
}

\fill (1.6, 1.3) circle (.05);
\draw (1.6, 1.3) circle (.5);

\end{scope}

\end{tikzpicture}
\caption{\label{fig:dscl} An example of a dynamic size cell list. The simulation box and cells are identically deformed. Larger cells are required compared to Figure~\ref{fig:docl} as a result of the shear.}
\end{figure}
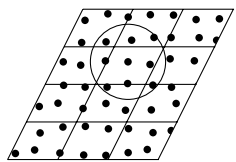

In the dynamic size cell list algorithm we choose the cells in the cell list to be congruent to the simulation box.  For convenience, the cells are all chosen to be of the same size, with the three heights of the cell all at least as large as the interaction range, see Figure~\ref{fig:dscl}.

\begin{figure}[ht]
\centering
\begin{tikzpicture}
\draw (0,0)--(1,0)--(1.5,1)--(.5,1)--(0,0);
\draw (.5,0)--(.5,1);
\draw(1,0)--(.2,.4);
\end{tikzpicture}
\caption{\label{fig:bh} Deformed simulation box with minimum box heights. When dividing the box into cells these heights determine the number of cells that can be placed in each direction.}
\end{figure}
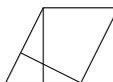

Let $c_i$ denote the number of cells in direction $i$, we can compute this given the box heights as follows:
\begin{align*}
c_i &= \left\lfloor \dfrac{h_i}{d_{cut}} \right\rfloor,
\end{align*}
where $h_i$ denotes the simulation box height in direction $i,$ see Figure~\ref{fig:bh}, and $d_{cut}$ denotes the cutoff distance for the interparticle interaction.  Once the simulation box heights and thus the number of cells in each direction have been determined the algorithm operates much as the equilibrium cell list algorithm, with each  cell neighborhood is composed of 27 cells.  It is convenient to convert particle positions to the underlying lattice basis of the simulation box when assigning particles to cells. The cell list is recomputed at each step using the newly deformed box. 
As the simulation box becomes more deformed the volume of each cell neighborhood grows as does the computational cost of computing pairwise interactions.  Bounds for the volume can be derived by examining the bounds on box deformation implied by the algorithm in use, Lees-Edwards, KR, or generalized KR.

\subsection{Dynamic Offset}
\label{subsec:do}

Choosing the cells to be congruent with the simulation box is natural, but not required. In the dynamic offset cell list algorithm, the cells maintain a rectangular shape. This is much like the Lees-Edwards sliding brick algorithm and faces some of the same problems \cite{Lees}.

\begin{figure}[ht]
\centering
\begin{tikzpicture}
\begin{scope}
\draw (-1,-1)--(1,-1)--(1.7,1)--(-.3,1)--(-1,-1);
\clip (-1,-1)--(1,-1)--(1.7,1)--(-.3,1)--(-1,-1);
\draw[step=.5 cm,black]
(-1,-1) grid (2,2);
\foreach \x in {-.9, -.5, -.1, .7, 1.1, 1.5, 1.9}
\foreach \y in {-.9, -.6, -.3, 0, .3, .6, .9}
{
\fill (\x + rand/20, \y + rand/20) circle (.05);
}

\foreach \x in {.3}
\foreach \y in {-.9, -.6, -.3, .3, .6, .9}
{
\fill (\x + rand/20, \y + rand/20) circle (.05);
}

\fill (.3, 0) circle (.05);
\draw (.3, 0) circle (.5);

\end{scope}
\end{tikzpicture}
\caption{\label{fig:docl} An example of a dynamic offset cell list. Here the cells do not deform with the simulation box. Note the cell heights are chosen so that the bottom vector of the box is neatly divided into cells as is the vertical box height.}
\end{figure}
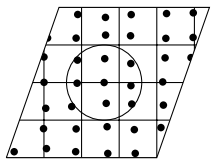

One problem introduced is dividing the deformed simulation box into rectangular cells. Allowing for partial cells along the faces would solve this problem, but these cells could be arbitrarily small and would not work well with periodic boundary conditions when choosing a cell neighborhood to ensure that all nonzero interactions are included. A preferable alternative is to combine partial cells along the edges into complete cells.

First we will briefly examine the 2D case. Here the deformed simulation box is a parallelogram. We assume without loss of generality that the parallelogram is oriented in space such that one corner lies at the origin and one side lies on the x axis and has length $l$. To find a division of this box into rectangular cells we can look at a rectangular rearrangement of the box, which we find by replacing the x coordinate of all points with $x\text{ mod }l$.

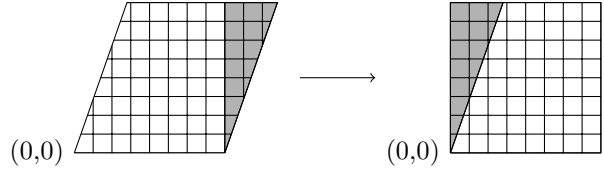
\begin{figure}[ht]
\centering
\begin{tikzpicture}
\begin{scope}[xshift = -2.5cm]
\filldraw[fill=black!30!white] (1,-1)--(1.7,1)--(1,1)--(1,-1);
\draw (-1,-1)--(1,-1)--(1.7,1)--(-.3,1)--(-1,-1);
\node at (-1.5,-1) {(0,0)};
\clip (-1,-1)--(1,-1)--(1.7,1)--(-.3,1)--(-1,-1);
\draw[step=.25 cm,black]
(-1,-1) grid (2,2);
\end{scope}

\begin{scope}[xshift = 2.5cm]
\filldraw[fill=black!30!white] (-1,-1)--(-.3,1)--(-1,1)--(-1,-1);
\draw (-1,-1)--(1,-1)--(1,1)--(-1,1)--(-1,-1);
\draw (-1,-1)--(-.3,1);
\draw[step=.25 cm,black]
(-1,-1) grid (1,1);
\node at (-1.5,-1) {(0,0)};
\end{scope}

\begin{scope}
\draw[->] (-.5, 0) -- (.5, 0);
\end{scope}
\end{tikzpicture}
\caption{\label{fig:rear}  We rearrange the simulation box to be rectangular by modding the x positions of all cells by the length of the x axis.}
\end{figure}

\begin{figure}[ht]
\centering
\begin{tikzpicture}
\begin{scope}
\draw [very thick](-1,-1)--(1,-1)--(1.6,1)--(-.4,1)--(-1,-1);
\clip (-1,-1)--(1,-1)--(1.6,1)--(-.4,1)--(-1,-1);
\draw[step=.25 cm,black]
(-1,-1) grid (2,2);
\filldraw [fill=black] (.25, .75) rectangle (.5, 1);
\filldraw [fill=black!30!white] (0, .5) rectangle (.25, .75);
\filldraw [fill=black!30!white] (0.25, .5) rectangle (.5, .75);
\filldraw [fill=black!30!white] (.50, .5) rectangle (.75, .75);
\filldraw [fill=black!30!white] (0, .75) rectangle (.25, 1);
\filldraw [fill=black!30!white] (.5, .75) rectangle (.75, 1);
\end{scope}
\begin{scope}[xshift = .6cm, yshift = 2cm]
\draw [very thick](-1,-1)--(1,-1)--(1.6,1)--(-.4,1)--(-1,-1);
\clip (-1,-1)--(1,-1)--(1.6,1)--(-.4,1)--(-1,-1);
\draw[step=.25 cm,black]
(-1,-1) grid (2,2);
\filldraw [fill=black!30!white] (.0, -1) rectangle (.25, -.75);
\filldraw [fill=black!30!white] (-.25, -1) rectangle (0, -.75);
\filldraw [fill=black!30!white] (-.5, -1) rectangle (-.25, -.75);
\filldraw [fill=black!30!white] (-.75, -1) rectangle (-.5, -.75);
\end{scope}
\begin{scope}[xshift = 2cm]
\draw [very thick](-1,-1)--(1,-1)--(1.6,1)--(-.4,1)--(-1,-1);
\clip (-1,-1)--(1,-1)--(1.6,1)--(-.4,1)--(-1,-1);
\draw[step=.25 cm,black]
(-1,-1) grid (2,2);
\end{scope}
\begin{scope}[xshift = 2.6cm, yshift = 2cm]
\draw [very thick](-1,-1)--(1,-1)--(1.6,1)--(-.4,1)--(-1,-1);
\clip (-1,-1)--(1,-1)--(1.6,1)--(-.4,1)--(-1,-1);
\draw[step=.25 cm,black]
(-1,-1) grid (2,2);
\end{scope}
\end{tikzpicture}
\caption{\label{fig:os} An example of an offset shift. The black cell is the cell under consideration and the surrounding gray cells together with the black cell form the cell neighborhood. Here as the cell was along the top edge of the simulation box, the top row of the neighborhood had to be expanded. This is because the cells in the simulation box were slightly offset to those in its periodic replica.  }
\end{figure}
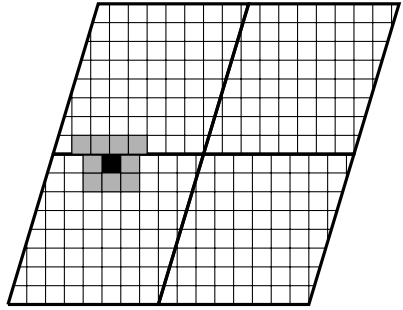

When assigning cells to a cell neighborhood in most cases nine cells are sufficient as in the equilibrium case. However cells which lie on the top or bottom side of the box may be fractionally offset from cells outside the box, requiring a larger cell neighborhood to ensure that all nonzero interactions are included.  In this case, ten cells will be needed, as seen in Figure~\ref{fig:os}. 

Moving to the 3D case we apply the same general method as above. We begin by finding the cubic division of the simulation box. We assume without loss of generality that the parallelepiped is rotated so one face is contained within the xy plane, one edge lies on the x axis, and one corner lies at the origin. This is equivalent to finding the QR decomposition of the simulation box where the shifted box lattice basis vectors are
\begin{equation*}
\label{eq:R}
R = \left[
\begin{array}{rrr}
r_{11} & r_{12} & r_{13} \\
0 & r_{22} & r_{23} \\
0 & 0 & r_{33} 
\end{array}
\right] 
= \left[
\begin{array}{rrr}
\vec{r_x} & \vec{r_y} & \vec{r_z}
\end{array}
\right].
\end{equation*} 
The rotated simulation box can be rearranged into a rectangular prism similarly to Figure~\ref{fig:rear}. For all points we reset their coordinates as follows
\begin{align}
y &= y- \left\lfloor\frac{y}{r_{22}}\right\rfloor r_{22} \\
x &=  \left( x- \left\lfloor\frac{y}{r_{22}}\right\rfloor r_{12} \right) \text{ mod }r_{22}
\end{align}
Using this method an interior cell neighborhood will be identical to one in the 3D equilibrium case, containing a total of 27 cells. As with the 2D case special handling is required for defining the neighborhood of cells lying on faces affected by offset shifts in order to ensure that all nonzero pairwise interactions are accounted for. In the 3D case we denote a face IJ if it is formed by vectors $\vec{r_i} \text{ and } \vec{r_j}$. Due to our reorientation of the simulation box no offset shifts are required for moving between replicas in the x direction.  Moving between replicas in the y direction only requires accounting for x directional offset shifts, as the box was oriented to lie flush with XY plane.  Moving between replicas in the z direction involves both x and y offset. Therefore, only the cells on the XZ or XY faces require extra work.

In handling the special cases we must expand the cell neighborhood as a result of imperfect cell alignment caused by periodic shifts. While this adds complexity to the algorithm, we show below that it results in an overall more efficient computation of particle interactions.

When determining the neighbor list for a cell there are four cases to consider:
\begin{enumerate}
\item{\bf{The cell does not lie on the XZ or XY faces of the simulation box:}} \\
In this case the cell neighborhood remains the same as in the equilibrium case as there are no unaccounted for shifts between periodic replicas.  The cell neighborhood has 27 cells.
\item{\bf{The cell lies on one of the XZ faces of the simulation box, but not one of the XY faces:}} \\
Here the cell neighborhood must be expanded to $4 \times 3$ cells on the XZ layer outside the simulation box as there could be a lattice shift in the Y direction resulting in imperfect cell alignment. A total of 30 cells are required for the cell neighborhood of these cells.
\item{\bf{The cell lies on one of the XY faces of the simulation box, but not one of the XZ faces:}} \\
The XY layer outside the simulation box could be shifted in both the X and Y direction. As a result that layer must be expanded to $4 \times 4$. A total of 34 cells are required for the cell neighborhood of these cells.
\item{\bf{The cell lies on both a XZ face and a XY face:}} \\
Here we must expand one of the XZ layers and one of the XY layers. We expand each affected layer as described above, but note that one cell in the expanded XZ layer will also lie in the expanded XY layer. In total 36 cells will be required to complete the cell neighborhood.
\end{enumerate}

\section{Efficiency Comparison}
\label{sec:eff}
If we assume an even density of particles throughout the simulation box the number of pair interactions computed is proportional to the volume of the area considered. Thus a comparison of the efficiencies of the algorithms can be made by a comparison of the average volume of the cell neighborhoods.

To calculate this we must find the volume of the individual cells and, for the dynamic offset algorithm, the average number of cells included in the neighbor list. It is important to note that the volume and thus efficiency of cell neighborhoods for both methods are time dependent. Particularly for the dynamic size cell list algorithm, the cell neighborhood volume fluctuates tremendously. To give an effective comparison between the two methods we will average their performance over a long period of time in a sample simulation.  As the deformation of the simulation box is bounded, this long term averaging gives us a very good idea of behavior for any type of run in which there is significant deformation. 

Since the search volume scales with the interaction radius, we choose to compare the methods in terms of their search efficiency. This is the ratio of non-zero particle interactions over the total number of particle pairs within a cell neighborhood. To calculate this figure given the geometry of the simulation we take a ratio of the volumes of a sphere with radius $d_{cut}$ and the cell neighborhood volume. In the equilibrium case if the simulation box sides are evenly divisible by $d_{cut}$ the cell neighborhood will be a cube with side length $3 d_{cut}$. Thus, the search efficiency at equilibrium is  
\begin{equation}\frac{\frac{4}{3}\pi d_{cut}^3}{27 d_{cut}^3} = \frac{\frac{4}{3} \pi}{27} \approx 0.1551.\end{equation}

\subsection{Dynamic Size Cell List Volume}\label{subsec:dscl}
Use of the dynamic size variation guarantees a 27 cell neighbor list. The cells are identical to each other and evenly divide the simulation box, thus their volume is the volume of the simulation box divided by the number of cells, or
\begin{equation} V_{DS} = \frac{27V_{box}}{c_1c_2c_3}.\end{equation} 
 
\subsection{Dynamic Offset Volume}
\label{subsec:DOV}
The number of cells in each direction (here called $l_i$ to avoid confusion) is easy to calculate from the rotated lattice basis $R$ in~\eqref{eq:R}.  We take the diagonal entries of this matrix and divide them by $d_{cut}$ rounding down to the closest integer. 
\begin{align*}
l_i &= \left\lfloor\frac{r_{ii}}{d_{cut}} \right\rfloor
\end{align*}
where $i$ is the direction considered and $r_{ii}$ is the appropriate entry in $R$.

As the cells are rectangular their volumes are the product of the side lengths. To find these we take the lengths of the simulation box used to size the cell list and divide them by the number of cells in the relevant direction as follows:
\begin{equation}V_{DO} = \frac{r_{11}r_{22}r_{33}}{l_1 l_2 l_3} \end{equation}

The number of cells in the average neighbor list requires additional computation.  To calculate it we take the sum of the number of cells in each of the four cases described in Subsection~\ref{subsec:do}, which we denote $S_i$ for case $i$,  multiplied by the number of cells comprising the neighbor list in each case, we then divide by the total number of cells:
\begin{align*}
avg &= \frac{27 S_{1}+30 S_{2}+34 S_{3}+36 S_{4}}{c_{1}c_{2}c_{3}} \\
&= \frac{27(l_{1} l_{2} l_{3}-(2l_1l_3-4l_1)-(2l_1l_2-4l_1)-4l_1)+30(2l_1l_3-4l_1)+34(2l_1l_2-4l_1)+36(4l_1)}{l_1l_2l_3}\\
&=\frac{27 l_{1} l_{2} l_{3}+14 l_{1} l_{2}+6 l_{1} l_{3}-4 l_{1}}{l_{1} l_{2} l_{3}}\\
&= 27+\frac{14}{l_{3}}+\frac{6}{l_{2}}-\frac{4}{l_{2}l_{3}}.	
\end{align*}
The average number varies depending on the size of the simulation. Intuitively we see the number of added cells approximates the ratio of surface area to volume of the simulation box and is inversely proportional to the number of cells.  

\subsection{General Comparisons}
For both dynamic offset and dynamic size algorithms, the average volume of a cell neighborhood depends on the time dependent lattice $L_t.$  To compare these two methods we can examine the average volume over a long run time. When using generalized KR boundary conditions a single long run with a fixed flow can provide good information about any run done with a general flow as set of all possible cell deformations  is independent of the flow matrix, and for almost every background flow the set of deformations is fully explored as time goes to infinity.

\begin{figure}[ht]
\centering
\includegraphics[width=9cm]{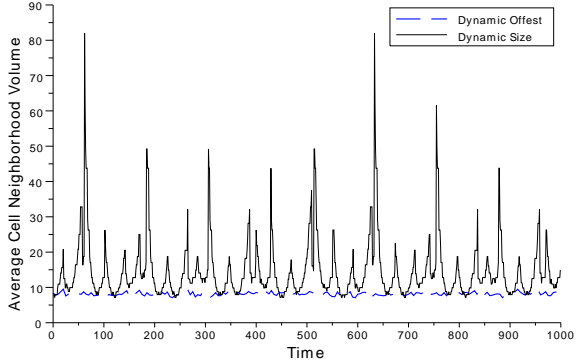}\\
\caption{Neighborhood size for a long run with a small simulation box with initial side lengths $a \approx 10 d_{cut}$. The black solid line is the cell neighborhood volume of the dynamic size cell list. The black dashed line is the neighborhood volume for the dynamic offset cell list.  The volumes have been normalized by choosing  $d_{cut} = \left( \frac{3}{4 \pi} \right)^{1/3}$ so the volume of a particle's interaction range is identically one.   Due to the large size of the simulation box the ratio of cells on the faces to interior cells is very small and thus the dynamic size cell list is quite efficient.}
\label{fig:eff}
\end{figure}

There is a complication in computing the search efficiency of the dynamic offset algorithm, the average number of cells in the cell list is dependent upon the overall size of the simulation box. We see from Subsection~\ref{subsec:DOV} that as the size of the simulation box grows this number goes to 27. The asymptotic search efficiency of the dynamic offset algorithm equals that of the equilibrium cell list algorithm, approximately $0.1551$, as the cell size and shape is determined identically to the equilibrium algorithm. In smaller simulations the search efficiency will be lower as the average number of cells in the cell list will be higher.   We can find the average volume of the cell list by averaging over a long run, which we perform with a uniaxial stretching flow,
\begin{equation}
A = \left[
\begin{array}{rrr}
0.05 & 0 & 0 \\
0 & -0.025 & 0 \\
0 & 0 & -0.025 
\end{array}
\right].
\end{equation}
In Figure~\ref{fig:eff}, we plot the average volume of the cell neighborhoods as a function of time for each algorithm.  The cutoff is scaled so the interaction range of a particle is one, $d_{cut} = \left( \frac{3}{4 \pi} \right)^{1/3}.$  Averaging the run time from Figure~\ref{fig:eff} gives us a search efficiency of approximately $0.0688$ for the dynamic size algorithm, compared to approximately $0.123$ for dynamic offset.
  
In comparisons of runtime efficiency the differences between the two methods applied to a small system were not quite as pronounced. In a run of 10,000 steps on a single core machine for 1728 particles with simple repulsive WCA interactions, we found the dynamic offset method took $71.7 \%$ the amount of time as the dynamic size method. Differences between the actual and idealized performance are partially explained by the size of the simulation used. Under the dynamic offset algorithm the simulation box was divided into about 10 cells in each direction.  Recomputing the estimated efficiency for this size we found it should have taken approximately $56 \%$ as long as the dynamic size method. The remaining difference can be explained as program overhead. Use of cell lists prevented unnecessary distance computations, but the number of force computations (which are more expensive) remained the same as force calculations were preceded with a distance check.
For large simulations, the additional overhead of the dynamic offset algorithm is balanced by the reduced number of particle interaction computations giving a more efficient algorithm overall.

\section*{Acknowledgements}
We would like to thank the generous support of Joan Barksdale, 
whose donations funded the REU project of the second and third authors.


\begin{thebibliography}{10}

\bibitem{bara99}
A.~Baranyai and P.~T. Cummings.
\newblock Steady state simulation of planar elongation flow by nonequilibrium
  molecular dynamics.
\newblock {\em J. Chem. Phys.}, 110(1):42--45, 1999.

\bibitem{dobs14}
M.~Dobson.
\newblock Boundary conditions for long-time nonequilibirum molecular dynamics
  simulations of incompressible flows.
\newblock {\em J. Chem. Phys.}, 18, 2014.
\newblock arXiv:1408.7078.

\bibitem{Evans}
D.~J. Evans and G.~P. Morriss.
\newblock {\em Statistical mechanics of nonequilibrium liquids}.
\newblock Anu E Press, 2007.

\bibitem{fren01}
D.~Frenkel and B.~Smit.
\newblock {\em Understanding Molecular Simulation}.
\newblock Academic Press, Inc., Orlando, FL, USA, 2nd edition, 2001.

\bibitem{hunt13}
T.~Hunt.
\newblock Periodic boundary conditions for the simulation of uniaxial
  extensional flow.
\newblock arXiv:1310.3905, 2013.

\bibitem{hunt10}
T.~A. Hunt, S.~Bernardi, and B.~D. Todd.
\newblock A new algorithm for extended nonequilibrium molecular dynamics
  simulations of mixed flow.
\newblock {\em J. Chem. Phys.}, 133(15):154116, 2010.

\bibitem{Kraynik}
A.~Kraynik and D.~Reinelt.
\newblock Extensional motions of spatially periodic lattices.
\newblock {\em In. J. Multiphase Flow}, 18(6):1045--1059, 1992.

\bibitem{Lees}
A.~Lees and S.~Edwards.
\newblock The computer study of transport processes under extreme conditions.
\newblock {\em J. Phys. C}, 5(15):1921, 1972.

\bibitem{Matin}
M.~Matin, P.~Daivis, and B.~Todd.
\newblock Cell neighbor list method for planar elongational flow: rheology of a
  diatomic fluid.
\newblock {\em Comput. Phys. Commun.}, 151(1):35--46, 2003.

\bibitem{todd99}
B.~Todd and P.~J. Daivis.
\newblock A new algorithm for unrestricted duration nonequilibrium molecular
  dynamics simulations of planar elongational flow.
\newblock {\em Comput. Phys. Commun.}, 117(3):191 -- 199, 1999.

\bibitem{todd98}
B.~D. Todd and P.~J. Daivis.
\newblock Nonequilibrium molecular dynamics simulations of planar elongational
  flow with spatially and temporally periodic boundary conditions.
\newblock {\em Phys. Rev. Lett.}, 81:1118--1121, Aug 1998.

\end{thebibliography}
\end{document}